\documentclass[11pt]{article}

\title{Uniform Versions of Infinitary Properties in 
 Banach Spaces}
\author{Carlos  Ortiz \\ Beaver College}

\usepackage{amsmath,amsthm,amsfonts}

\chardef\bslash=`\\ 

\newcommand{\findef}{\mbox{\rule{2.0mm}{2.0mm} }}

 \theoremstyle{plain} 
\newtheorem{thm}{Theorem}[section]
\newtheorem{cor}[thm]{Corollary}
\newtheorem{lem}[thm]{Lemma}

\theoremstyle{definition}
\newtheorem{defn}[thm]{Definition}
\newtheorem{example}[thm]{Example}

\theoremstyle{remark}
\newtheorem{rem}[thm]{Remark}
\newtheorem*{notation}{Notation}

\newcommand{\apptruth}{\mbox{$\models_{AP}$\/} }

\newcommand{\structa}{\mbox{\boldmath $E $\/} }

\newcommand{\structn}{\mbox{\boldmath $E_{n}$\/}}

\newcommand{\verify}{\mbox{$\models$ }}

\newcommand{\eval}[2][\right]{\relax
  \ifx\#1\right\relax \left.\fi\#2\#1\rvert}

\begin{document}
\maketitle

\begin{abstract}
In functional analysis it is of interest to study the
 following general question:

\begin{quote} Is the
uniform version of a property that holds  
in all Banach spaces also valid in all Banach spaces?
\end{quote}

Examples of affirmative answers to the above question are the host of
 proofs of
almost-isometric versions of well known isometric theorems. Another
example is  
Rosenthal's
uniform version of Krivine's Theorem. Using
an extended version of
 Henson's Compactness result for positive bounded formulas in normed
structures, we show
that the answer of the above question is in fact yes for 
 every property that can be expressed in a particular infinitary
 language. Examples of applications are given.
\end{abstract}

\section{Introduction}

A natural type of questions in functional analysis asks if the
 ``almost'' version of a
theorem true in a class of normed spaces is also true in the class.
Here are some examples of such questions:

\begin{enumerate}
\item Ulam's Theorem:
\begin{itemize}
\item
\begin{thm}Let $T:X\rightarrow Y$ be an onto 
function from a Banach space $X$ to a Banach space $Y$
with $T(0)=0$ such that:
$$||T(x)-T(y)|| =||x-y|| \mbox{ for all } x, y \in X$$
then $||T(x+y)-T(x)-T(y)||=0 \mbox{ for } x,y \in X$
\end{thm}
\item   
Gevirtz (\cite{Gevirtz}) proved the following "almost" version:
\begin{thm} Let $T:X\rightarrow Y$ be an onto function from a Banach
space $X$ to a Banach space $Y$ with $T(0)=0$ such that $\forall x,y
\in X$:
$$ (1-\epsilon) ||x-y|| \leq ||T(x)-T(y)|| \leq (1+\epsilon) ||x-y||$$
 then $||T(x+y)-T(x)-T(y)|| \leq \epsilon'(||x||+||y||) \mbox{ for } x, y \in X$,
where $\epsilon' \rightarrow 0$ as $\epsilon \rightarrow 0$.
\end{thm}
\end{itemize}
\item A classical result of Behrends (\cite{Behrends}):
\begin{itemize}
\item A linear projection $P:E\rightarrow E$ is called an
 $L^{p}$-projection, 
$1 \leq p \leq \infty$, if $\forall x \in E$,
$$ (||P(x)||^{p} + ||x-P(x)||^{p})^{1/p}=||x||,$$ with the obvious
modification for the case $p=\infty$.

Behrends proved the following:
\begin{thm}
Let $E$ be a Banach space with dim ($E$) $> 2$. Let $ 1 \leq p,q \leq \infty$,
$p \neq 2$ and  such that $P,Q:E\rightarrow E$ are $L^{p}$ and $L^{q}$
projections. Then $p=q$ and $||PQ-QP||=0$. 
\end{thm}
\item
The "almost" isometric case was proved by Cambern, Jaroz and Wodinski
 (\cite{Cambern}):
\begin{thm}
Let $E$ be a Banach space with dim $(E) > 2$. Let $1 \leq
p,q \leq \infty$, $p \neq 2$ and let $P,Q:E\rightarrow E$ be projections
with the additional properties:
$$\forall x \in E \mbox{, } (1-\epsilon)||x|| \leq ||P(x)||^{p}
+ ||x-P(x)||^{p})^{1/p} \leq (1+\epsilon)||x||$$
and
$$\forall x \in E \mbox{, } (1-\epsilon)||x|| \leq ||Q(x)||^{q}
+ ||x-Q(x)||^{q})^{1/q} \leq (1+\epsilon)||x||$$
then
$$|p-q|\leq \epsilon'(p) \mbox{ and } ||PQ-QP||\leq \epsilon'(p) ,
\mbox{ where } \epsilon' \rightarrow 0 \mbox{ as } \epsilon \rightarrow
0.$$
\end{thm}
\end{itemize}
 \end{enumerate}

Jarosz (\cite{Jarosz}) pointed out that the proof of the two
above results could be 
simplified considerably by using ultraproducts of Banach spaces. 
One may ask then  if it is possible to study this phenomena
 in a systematic way from a logical point of view.

The natural model theoretic setting to answer this question
is Henson's  logic
of positive bounded formulas in normed spaces (\cite{Henson}).
This logic
$L_{PB}$  is closed under
 finite conjunction, finite disjunction and
bounded quantification. The normed spaces that
 are the natural models for this language are called normed space
 structures. Henson defined the notion of an $n$-approximation
of a formula $\phi$ in $L_{PB}$, denoted by $(\phi)_{n}$. From this
concept he defined the semantic notion of
 approximate truth ($\apptruth$) for this logic. It can be seen that
 $(L_{PB},\apptruth)$ has a compactness theorem 
(see \cite{Henson} and \cite{Henson&Iovino}).

However, $L_{PB}$ has a fundamental limitation for our purposes:
formulas of the form
 $\phi \Rightarrow \psi$ (like Ulam's Theorem or Behrend's
result), with 
$\phi, \psi$ in $L_{PB}$, are not in $L_{PB}$.
Furthermore, since most interesting statements in functional analysis
are of a fully infinitary type (using countable disjunction for example),
 we are interested in obtaining a
 general uniformity result that includes
 infinitary formulas with countable disjunctions and
 conjunctions and bounded 
quantification over infinitely many variables.

We deal with this limitation by extending the notion of approximate truth
($\apptruth$) 
to an infinitary logic 
$L_{A}$ that contains $L_{PB}$ and is closed under countable
 conjunctions ($\bigwedge$), 
negation ($\neg$) and bounded existential quantification over
 countably many 
variables ( $\exists \vec x (\bigwedge_{i=1}^{\infty} ||x_{i}||
\leq r_{i}\wedge \ldots) $). This logic   and
 its corresponding notion of $\apptruth$ was introduced in
 \cite {Ortiz} to study the idea of
``proof by approximation'' in analysis from a logical point of view. 

Let us describe briefly how we extend the notion of
approximation of a formula from $L_{PB}$ to $L_{A}$.
The main obstacle to the extension of  approximate truth to $L_{A}$
 is the negation connective. The key  to solving this
is to 
extend  Henson's idea   
 of a sequence of approximate formulas $\{(\phi)_{n}:n \in \omega\}$ 
 (for formulas $\phi \in L_{PB}$) to a tree of  positive bounded formulas 
 $\{([\phi]_{h})_{n}: h \in I(\phi) \mbox{, } n \in \omega\}$ where each 
 sequence $([\phi]_{h})_{1},([\phi]_{h})_{2},\ldots,([\phi]_{h})_{n},
\ldots$ is a branch
 of  the
 tree. In this way a
 sentence $\phi$ is 
 approximately true in a normed structure $E$ iff there exists a branch 
 of the tree of approximations of $\phi$ such that all the approximations 
 of this branch hold in $E$. An analogue to this approach in
 classical infinitary logic is the notion of approximation of infinitary
formulas by Vaught sentences (
 \cite{Hodges}).

Since the approximations for the infinitary formulas in $L_{A}$ are
 in turn positive bounded formulas, we
can invoke Henson's theorem for $L_{PB}$
to obtain a compactness result for $(L_{A}, \apptruth)$.
From this  compactness theorem we get the following general uniformity
 result for  formulas of the form $\neg \phi$:

\begin{quote} Uniformity Theorem for $L_{A}$.

For any class of normed space structures axiomatized by a
 theory $\Sigma$ in $L_{PB}$, 
for any sentence $\phi \in L_{A}$, if $\Sigma
 \models \neg \phi$ then for every branch  $h \in 
I(\phi)$ there exists an integer $n$ such that $\Sigma
 \models \neg ([\phi]_{h})_{n}$.
\end{quote}

This paper is organized as follows:
For the sake of completeness, the first two sections are
devoted to a brief review of
 Henson's notion
 of a normed space structure as well as the definition of the logic
$L_{A}$ (Section~\ref{LA}) and a review of the definition of
 $\apptruth$  in
$L_{A}$ (Section~\ref{apptruth}).

In Section~\ref{comply} we use Henson's Compactness theorem for 
$L_{PB}$ to get
a Model Existence Theorem for $(L_{A},\apptruth)$.
 From this theorem
we prove the Uniformity Theorem.

Finally, in Section~\ref{Applications} we give applications of the
 Uniformity Theorem  to Banach space theory. Among 
the applications we can cite the theorems of Gervitz and Jaroz
 mentioned above, as well
as Rosenthal's uniform version of 
Krivine's Theorem.

A note on notation: we will use \findef  to denote the end of
 definitions,
 examples and remarks.
 
\section{Normed space structures and the logic $L_{A}$  }
\label{LA}

We begin by briefly recalling  Henson's notions of normed space structure
 and of a language for normed space structures. For a more detailed account,
the reader may look at \cite{Henson} or \cite{Iovino}.

\begin{defn} Real valued m-ary relations.

A real valued m-ary relation on a normed space $E$ is a
 function $\cal R$$:
E^{m} \rightarrow \mathbb{R}$ which is uniformly continuous
 on every bounded
subset of $E^{m}$. \findef
\end{defn}
 
\begin{defn} A normed space structure is a structure of the form
\begin{equation*}
\structa=(E, f_{i},\mbox{$\cal R$}_{j}: i \in I, j \in J)
\end{equation*}
with:
\begin{itemize}
\item $E$ being a normed space structure over the reals;
\item each $f_{i}$ being a function $f_{i}:E^{m} \rightarrow E$
 for some natural number $m$;
\item every $f_{i}$ being a uniformly continuous function on every
 bounded subset of $E^{m}$;
\item each $\cal R$$_{j}$ being a real valued relation. \findef
\end{itemize}

\end{defn}

From the concept of normed space structure follows the notion of a
  signature  for normed space structures. 
A signature in this setting consists of function symbols, real valued
relation symbols,  bounds
 for the function symbols and the real valued
relation symbols. It also contains moduli of uniform continuity
 for the function and
real valued relation symbols on each bounded set.

\begin{notation} We will use  $||.||$ for norms, and  $|.|$
 for absolute value.
\end{notation}

\begin{defn}

A signature for the normed space structure
\begin{equation*}
\structa=(E, f_{i} ,\mbox{$\cal R$}_{j}: i \in I, j \in J)
\end{equation*}
 consists of:
 \begin{enumerate}
 \item an $m$-ary function symbol for each m-ary function $f_{i}$;
 \item an $m$-ary relation symbol for each $m$-ary real valued
 relation $\cal R$$_{i}$;
\item for each $m$-ary function $f_{i}$, each
positive integer $N$ and each positive rational $\epsilon$,
 a positive rational
$\delta(f_{i}, N,\epsilon)$ such that
\begin{equation*}
\begin{split}
||x_{k}||, ||y_{k}|| < N \mbox{ and } ||x_{k}-y_{k}||
 < \delta(f_{i}, N, \epsilon)
\mbox{ } (1 \leq k \leq m) \\
\mbox{implies } ||f_{i}(x_{1},\ldots,x_{m}) -
 f_{i}(y_{1},\ldots,y_{m})|| < \epsilon \mbox{;}
\end{split}   
\end{equation*}
\item for each $m$-ary real valued function $\mbox{$\cal R$}_{j}$, each
positive integer $N$ and each $V \in \cal V$, a positive rational
$\delta(\mbox{$\cal R$}_{j}, N,V)$ such that
\begin{equation*}
\begin{split}
||x_{k}||, ||y_{k}|| < N \mbox{ and } ||x_{k}-y_{k}||
 < \delta(\mbox{ $\cal R$}_{j}, N, V)
\mbox{ } (1 \leq k \leq m) \\
\mbox{ implies } (\mbox{$\cal R$} _{j}(x_{1},\ldots,x_{m}) -
\mbox{$\cal R$} _{j}(y_{1},\ldots,y_{m})) \in V \mbox{;}
\end{split}   
\end{equation*}
\item for each $m$-ary function $f_{i}$, for each integer $N$, an integer $K(i,N)$ such that 
\begin{equation*}
||x_{k}|| \leq N \mbox{ } (1 \leq k \leq m) \mbox{ implies } ||f(\vec x)|| \leq K(i,N);
\end{equation*}
\item for each $m$-ary relation $\cal R$$_{j}$, for each integer $N$, an integer 
$K(j,N)$ such that:
\begin{equation*}
||x_{k}|| \leq N \mbox{ }(1\leq k\leq m) \mbox{ implies }
 |\mbox{$\cal R$}_{j}(\vec x)| \leq K(i,N).\mbox{ } \findef
\end{equation*}
\end{enumerate}
\end{defn}

If $\Omega$ is a signature for the normed space structure
 $\structa$, we say that $\structa$ is an $\Omega$-structure.

\begin{rem}
A constant in this signature is a $0$-ary function. \findef
\end{rem}
 
For every signature for normed structures, a first order
 language is associated in the following way.

\begin{defn}
Let $\Omega$ be a signature for a normed space structure with universe $E$. We associate
to it the following first order language consisting of:
\begin{itemize}
\item a constant symbol $0$, a binary function symbol $+$, and for each rational
 scalar $r$, a function symbol for the scalar multiplication $x\rightarrow rx$;
 \item for each rational number $r$, predicate symbols for the sets
 \begin{equation*}
 \{x \in E: ||x || \leq r\} \mbox{ and } \{x \in E: ||x|| \geq r\}
 \end{equation*}
\item the function symbols of $\Omega$;
\item for each real valued  relation symbol $\cal R$ in $\Omega$, for each
 rational $r$, 
predicate symbols for the sets:
\begin{equation*}
\begin{split}
\{(x_{1},\ldots,x_{m}) \in E^{m}: \mbox{$\cal R$}(x_{1},\ldots,x_{m})
 \leq r \} \mbox{ and } \\
 \{(x_{1},\ldots,x_{m}) \in E^{m}:
\mbox{$\cal R$}(x_{1},\ldots,x_{m})
 \geq r \}
\mbox{. \findef}
\end{split}
\end{equation*}
\end{itemize}

\end{defn}
The formulas of $L_{PB}$ are defined by induction.
As usual, for every formula $\phi$, 
we will use the notation $\phi(\vec x)$ to 
express the fact that the free variables of $\phi$ are among the components
of the vector $\vec x$. Likewise, $\phi(\vec x_{1},\vec
x_{2},\ldots)$ means that the free variables of $\phi$ are among the
components of the vectors $\vec x_{1},\vec x_{2},\ldots$.

\begin{defn} Definition of $L_{PB}$.

Fix a signature $\Omega$.
\begin{enumerate}
\item If $t$ is a term of the  first order language corresponding to
 $\Omega$ and
$r$ is a rational number, then $||t|| \leq r$ and $||t || \geq r$
 are formulas in $L_{PB}$.
\item Let $\cal R$ be an $m$-ary real valued relation,
 $t_{1},\ldots, t_{m}$
 be terms and $r$ be a rational number, then
 $\cal R$$(t_{1},\ldots,t_{m})
 \leq r$ and $\cal R$$(t_{1},\ldots,t_{m}) \geq r$ are
 formulas in $L_{A}$.

\item If $\phi_1,\phi_2$ are  formulas in $L_{PB}$, then $\phi_{1}
\wedge \phi_{2}$ is a formula in $L_{PB}$.
\item if $\{\phi_{i}\}_{i=1}^{\infty}$ is a countable collection of
formulas in $L_{PB}$, then $\bigwedge_{i=1}^{\infty} \phi_{i} \in
L_{PB}$.

\item If $\phi_{1}$, $\phi_{2}$ are formulas in $L_{PB}$ then
$\phi_{1} \vee \phi_{2}$ is a formula in $L_{PB}$. 

\item Consider a  formula $\phi( y,
\vec x)$ in $L_{PB}$. Let $r \geq 0 \in \mathbb{Q}$.
 The following formula is in $L_{PB}$: $
    \exists y( ||y || \leq r \wedge 
      \phi( y,\vec x ))$.

\item
Consider a  formula $\phi( y.
\vec x)$ in $L_{PB}$. Let $ r \geq 0 \in \mathbb{Q}$.
 The following formula is in $L_{PB}$:
  $
    \forall y
 ( ||y || \leq r \Rightarrow 
      \phi( y,\vec x ))$ \findef
\end{enumerate}
\end{defn}

{\bf Note:\/} For real  valued relations $\cal R$$_{1}(\vec x)$ and $\cal R$$_{2}(\vec x)$, 
we will write $\cal R$$_{1}(\vec x) \leq \cal R$$_{2}(\vec x)$ in $L_{PB}$ to abbreviate 
the formula  $\bigwedge_{q\in \mathbb{Q}} \cal R$ $_{1}(\vec x) \leq q \vee \cal R$ $_{2}(\vec x) \geq q$. In a similar manner we will abbreviate 
$\cal R$$_{1}(\vec x) \geq \cal R$$_{2}(\vec x)$ and $\cal R$$_{1}(\vec x) = \cal R$$_{2}(\vec x)$.

We now recall the definition of approximate formulas for
$L_{PB}$ (see \cite{Henson} for more details).

\begin{defn} Definition of $\apptruth$ for $L_{PB}$.
\label{logichenson}

 For every $\phi \in L_{PB}$ and every integer $n$, define
 $\phi_{n}$ as follows:
\begin{itemize}
\item 
$(||t|| \leq r)_{n}: ||t|| \leq r + 1/n$ and
 $(|| t||  \geq  r)_{ n} :  ||t|| \geq  r - 1/n$.
Likewise, $(\cal R$$(\vec t) \leq r)_{n} : \cal R$$(\vec t)
\leq r+1/n$ and 
$(\cal R$$(\vec t) \geq r)_{n} : \cal R$$(\vec t) \geq r-1/n$.
\item $(\phi \wedge \psi)_{n}: \phi_{n} \wedge \psi_{n}$ and 
 $(\phi \vee \psi)_{n}: \phi_{n} \vee \psi_{n}$.
 \item $(\bigwedge_{i=1}^{\infty} \phi_{i})_{n} :
 \bigwedge_{i=1}^{n} (\phi_{i})_{n}$.
 \item $( \exists y (||y|| \leq r \wedge
\phi(\vec x, y)))_{n} :
 \exists y (||y|| \leq  r+1/n \wedge \phi_{n}(\vec x, y))$.
 \item $(\forall y(||y|| \leq \vec r \Rightarrow
\phi(\vec x,y)))_{n}
 : \forall y(|| y|| \leq  r-1/n \Rightarrow
\phi_{n}(\vec x,y))$.
 \end{itemize}

Finally, for $\phi \in L_{PB}$, $\structa \apptruth
 \phi(\vec a)$
iff $\structa \verify \bigwedge_{n=1}^{\infty} \phi_{n}(\vec a)$.\findef
\end{defn}

{\bf Notation:\/} to avoid long formulas, we will abbreviate
$$\exists x_{1}(||x_{1}|| \leq r_{1} \wedge \exists x_{2} (
||x_{2}|| \leq r_{2} \wedge \ldots \exists x_{s}(||x_{s} ||
\leq r_{s} \wedge \phi)\ldots )$$
by $\exists \vec x(\bigwedge_{i=1}^{s}
 ||x_{i} || \leq r_{i} \wedge \phi)$.

It is easy to see (\cite{Henson}) that the following is true:
\begin{thm} 
\label{natural}

Let $E$ be a normed space structure, and let
 $\phi(\vec x)$ be a positive bounded 
formula. The following holds for every
normed structure $E$:
\begin{itemize}
\item For every integer $n$, $E \models \phi_{n+1}(\vec a)
\Rightarrow \phi_{n}(\vec a)$;
\item If  $E \models \phi(\vec a) \mbox{ then }
E \apptruth \phi(\vec a)$.
\end{itemize}
\end{thm}

Although $L_{PB}$ does not have the negation connective, one
 can define a weak approximate
negation operator in
$L_{PB}$ inspired by Henson's weak negation operator
(\cite{Henson&Iovino}).

\begin{defn} Weak approximate negation operator.

Fix a signature $\Omega$. For every integer $n$ and every
formula $\phi \in L_{PB}$ we define the operator
$neg(\phi,n)$ as follows:
\begin{enumerate}
\item If $t$ is a term of the  first order language corresponding to
 $\Omega$ and
$r$ is a rational number, then $neg (||t|| \leq r, n):
||t|| \geq r+1/n$ and $neg (||t || \geq r,n): ||t|| \leq r-(1/n)$.
Likewise, $neg(\cal {R}$$(\vec t)\leq r,n): \cal {R}$$(\vec t)
\geq r+1/n$ and $neg(\cal {R}$$(\vec t)\geq r,n): \cal {R}$$(\vec t)
\leq r-1/n$.
\item $neg (\phi\wedge \psi, n): neg(\phi,n) \vee neg (\psi,n)$.
\item $neg (\phi\vee \psi, n): neg(\phi,n) \wedge neg (\psi,n)$.
\item $neg(\bigwedge_{i=1}^{\infty} \phi_{i}, n) : \bigvee_{i=1}^{n}
 neg (\phi_{i},n)$.
\item $neg(\exists x(|| x|| \leq r
\wedge \phi,n) : \forall  x (|| x|| \leq  r + (1/n)
\Rightarrow neg(\phi,n))$.
\item $neg(\forall x(||x|| \leq  r
\Rightarrow \phi,n) : \exists  x (||x|| \leq \vec r - (1/n)
\wedge neg(\phi,n))$.\findef
\end{enumerate}
\end{defn}

The main property of the weak approximate negation operator
 is given by the
following lemma. We call the subcollection of $L_{PB}$ containing the
atomic formulas and closed under finite conjunction, disjunction, and
the existential and universal bounded quantification the {\bf finitary
part of $L_{PB}$\/}.

\begin{lem} \label{negando}
For every formula $\phi \in L_{PB}$,
\begin{itemize}
\item $\forall n \in \omega$, $neg (\phi,n)$ is in the finitary
part of $L_{PB}$.
\item For every structure $E$ and every $\vec a \in E$,
$$ E \not \apptruth \phi(\vec a) \mbox{ iff } \exists m \in \omega
\mbox{ } E \models neg(\phi,m).$$
\item For every integer $n$, for every structure $E$ and every
$\vec a $ in $E$, $$E \models \neg (\phi_{n}) \Rightarrow neg(\phi,n+1).$$
\end{itemize}
\end{lem}

\begin{proof} The proof is direct and is left to the reader. \end{proof}

We now define the fully infinitary logic $L_{A}$ based on $L_{PB}$.

 \begin{defn} Definition of $L_{A}$.

Fix a signature $\Omega$. We define $L_{A}$ by induction in
formulas:
\begin{enumerate}
\item $L_{PB} \subset L_{A}$.

\item If $\phi_1,\phi_2,\ldots,\phi_i$\ldots ($i{<}\omega$)
is a collection of formulas in $L_{A}$,
then for every integer 
$n$,  $\bigwedge_{i=1}^{n} \phi_{i} \mbox{ , and }
 \bigwedge_{i=1}^{\infty}
\phi_{i} $ 
are formulas in $L_{A}$.
\item If $\phi$ is a formula in $L_{A}$ then $\neg\phi$ is
also a formula in $L_{A}$.

\item Consider a  formula $\phi( y_{1},\ldots , y_{n},\ldots,
\vec x)$ in $L_{A}$. Let $\vec r=(r_{1},\ldots, r_{n}, \ldots)$ be a
 corresponding vector
 of rational numbers.
 The following formula is in $L_{A}$:
  \begin{displaymath}
    \exists (y_{1},\ldots , y_{n},\ldots )
 (\bigwedge_{n=1}^{\infty} ||y_{n} || \leq r_{n} \wedge 
      \phi( y_{1},\ldots , y_{n},\ldots ,\vec x )). \mbox{ \findef }
        \end{displaymath}
\end{enumerate}
\end{defn}

{\bf Notation:\/} to avoid very long formulas we will abbreviate
\begin{displaymath}
\exists ( y_{1},\ldots,  y_{n},\ldots ) 
(\bigwedge_{n=1}^{\infty} ||y_{n}|| \leq r_{n}  \wedge
\phi(\vec x ,\vec y))
\end{displaymath}
by $\exists \vec y(||\vec y|| \leq \vec r \wedge \phi(\vec x,\vec y))$.
 Likewise $\neg
\exists \vec y (||\vec y || \leq \vec r \wedge  \phi(\vec y,\vec x))$
 will be
abbreviated by $\forall \vec y(||\vec y || \leq \vec r  \Rightarrow
 \neg \phi(\vec
y,\vec x))$. We will also abbreviate $\neg \bigwedge \neg$ by
 $\bigvee$ and $\neg
(\phi \wedge \neg \psi)$ by $ \phi \Rightarrow \psi$.

Finally, given a countable set $A=\{a_{1},\ldots,a_{n},\ldots\}$ 
 with a fixed enumeration 
 and  countable formulas $\{\phi_{a}\}_{a\in A}$ we understand 
 by $\bigwedge_{a\in A} \phi_{a}$ the formula 
 $\bigwedge_{n=1}^{\infty} \phi_{a_{n}}$. 
 Likewise, for an  arbitrary integer $m$, 
 we understand by $\bigwedge_{a \in A\uparrow m} 
 \phi_{a}$ the formula $\bigwedge_{n=1}^{m} \phi_{a_{n}}$.

The notion of satisfaction ($\structa \models \phi(\vec a)$) for
$\Omega$-structures $\structa$,  with $\vec a$ a vector of
 elements in $E$ and for $\phi \in L_{A}$
is the natural one and we are not going to do it here.
 The interested reader is directed to \cite{Henson&Iovino}
 for more details.

\begin{example} Expressive Power of $L_{A}$.
\label{reflexive1}

We show that the property of reflexivity can be expressed in the
 logic $L_{A}$.

For any Banach space $(X,||.||)$, let $B_{1}$ denote the unitary ball.
 A characterization of reflexivity due to James (\cite{James})
(see also \cite{Van Dulst}) 
 that does not require any 
mention of the dual is the following:

{\em 
A Banach space $(X,||.||)$  is reflexive 
iff\/} 
\begin{equation*}
 \forall \epsilon > 0 \mbox{ }
\forall \{x_{i}\}_{i=1}^{\infty} \subseteq  B_{1}\mbox{ } \exists k
 \in \omega \mbox{ }
dist[conv(\{x_{1},\ldots,x_{k}\}),conv(\{x_{k+1},\ldots\})] 
\leq \epsilon
\end{equation*}
Here, for any set  $A \subseteq X$,  $conv(A)$ is the convex hull spawned 
by $A$. Similarly,  given two sets $A,B \in X$, $dist[A,B]=inf\{||x-y||
 : x \in A, y \in B\}$.

 Let $\Omega$ be the empty signature. Then the $\Omega$-structures are the normed
 spaces.
 The following sentence of $L_{A}$ expresses reflexivity for the
 closure of  these structures.
\begin{equation*}
\bigwedge_{n=1}^{\infty} \forall \vec x (|| \vec x || \leq 1 \Rightarrow
\bigvee_{k=1}^{\infty} \bigvee_{r=1}^{\infty} 
\bigvee_{\vec a \in CO(k)} 
\bigvee_{\vec b \in CO(r)} 
||\sum_{i=1}^{k} a_{i}x_{i} -\sum_{j=1}^{r} b_{j}x_{k+j}|| 
< (1/n)
\end{equation*}

Here, $\forall s \in \omega$, $CO(s)$ 
is the subset of $\mathbb{Q}^{s}$ made of all the $s$-tuples 
$(a_{1},\ldots , a_{s})$ such that $\sum_{i=1}^{s} a_{i}=1$ and 
$a_{1},\ldots , a_{s} \geq 0$. \findef
\end{example}

\section{Approximate formulas for $L_{A}$}
\label{apptruth}

Our intention is to generate approximations of all the formulas
in $L_{A}$ by using the  formulas in $L_{PB}$ as building blocks.
 As mentioned in the introduction, the main 
problem arises from the negation connective.
 We will use the weak approximate negation operator 
($neg(.,.,.)$) defined in the previous section to solve this problem.

Formally, we associate to every formula $\phi$ in $L_{A}$ a set of
indices $I(\phi)$ (the branches of the tree of approximate
formulas) and for every $h \in I(\phi)$ a
formula
$[\phi]_{h} \in L_{PB}$. Intuitively, for  every branch $ h\in I(\phi)$,
the approximate formulas of $[\phi]_{h}$ (the collection 
$\{([\phi]_{h})_{n}| n \in \omega\}$) are going to
``approach'' $\phi$ as $n$
tends to $\infty$. 

The notions of $I(\phi)$ and $([\phi]_{h})_{n}$ were introduced
 (in a different presentation) 
 in \cite{Ortiz}.

{\bf Notation\/}:  Given two formulas $\phi, \sigma$, we 
will write $\phi \equiv \sigma$ if $\phi$ and $\sigma$ are identical
 formulas.

In the rest of this section we fix a signature $\Omega$.

\begin{defn} Approximate formulas in $L_{A}$.
\label{defofapprox}

For any formula $\phi(\vec x)$ in $L_{A}$ we define by induction
 in formulas:
\begin{itemize}
\item a set $I(\phi)$ of branches;
\item $\forall h \in I(\phi)$, formulas $[\phi]_{h} \in L_{PB}$.
\end{itemize}
 {\bf Formulas in $L_{PB}$.\/} 
$\forall \phi \in L_{PB}$, $I(\phi)=\{\emptyset\}$. Furthermore,
$[\phi]_{\emptyset}: \phi$.

{\bf Countable (Finite) Conjunction.\/} For any countable (or finite)
collection 
$\{\phi_{i}\}_{i=1}^{\infty}$ ($\{\phi_{i}\}_{i=1}^{m}$)
of
formulas in $L_{A}$, we define:
\begin{itemize}
\item $I(\bigwedge_{i=1}^{\infty} \phi_{i}(\vec
x))=\prod_{i=1}^{\infty}
I(\phi_{i})$ (the cartesian product of the $I(\phi_{i})$)  (or
$I(\bigwedge_{i=1}^{m} \phi_{i}(\vec
x))=\prod_{i=1}^{m}
I(\phi_{i})$).

\item For every $h$ in $I(\bigwedge_{i=1}^{\infty} \phi_{i})$,
\begin{displaymath}
[\bigwedge_{i=1}^{\infty} \phi_{i}]_{h}:
\bigwedge_{i=1}^{\infty} [\phi_{i}]_{h(i)} \mbox{ ( or
} [\bigwedge_{i=1}^{m} \phi_{i}]_{h}:
\bigwedge_{i=1}^{m} [\phi_{i}]_{h(i)})
\end{displaymath}
\end{itemize}

{\bf Negation.\/} For any formula $\phi$ in $L_{A}$, we have:
\begin{itemize}
\item $I(\neg \phi) \subseteq (I(\phi)\times \omega)^{\omega}$
is the collection of all maps $f=(f_{1},f_{2})$ with the
following ``weak'' surjectivity property:
\begin{displaymath}
\forall h \in I(\phi) \mbox{ } \exists s \in \omega \mbox{, }
([\phi]_{h})_{f_{2}(s)} \equiv ([\phi]_{f_{1}(s)})_{f_{2}(s)}
\end{displaymath}
\item For every $f=(f_{1},f_{2}) \in I(\neg \phi)$, $ [\neg \phi]_{f} : \mbox{ } \bigwedge_{s=1}^{\infty} 
neg([\phi]_{f_{1}(s)},f_{2}(s))$
\end{itemize}

{\bf Existential.\/} For every formula $\phi(\vec y,\vec
x)$, for every corresponding vector $\vec r$ of rational numbers,
 we have:
\begin{itemize}
\item $I(\mbox{ }\exists \vec y(||\vec y|| \leq \vec r \wedge
\phi(\vec y, \vec x)\mbox{)
})=I(\phi(\vec y, \vec x))$.
\item For every $h$ in $I(\mbox{ }\exists \vec y(||\vec y || \leq \vec r \wedge
 \phi(\vec y,\vec x))\mbox{ })$, let $Ind(n)$ be the value of the maximal
index such that $x_{In}$ appears free in
 $([\phi(\vec y, \vec x)]_{h})_{n}$.
We define
\begin{equation*}
[\exists \vec y (||\vec y || \leq \vec r \wedge
\phi(\vec y, \vec x))\mbox{ }]_{h} : \bigwedge_{n=1}^{\infty}
\exists \vec y (\bigwedge_{s=1}^{Ind(n)} ||y_{i}|| \leq  r_{i} \wedge
([\phi(\vec y,\vec x)]_{h})_{n} \mbox{ }) 
\mbox{. \findef  } 
\end{equation*}
\end{itemize}

\end{defn}

 The
formulas $([\phi]_{h})_{n}$ are the {\bf approximate formulas\/}  of
$\phi$. 

\begin{defn} Approximate Truth.

Fix an $\Omega$-structure $\structa$. Let $\phi(\vec x)$ be an arbitrary 
formula in $L_{A}$. 
 We say that $
 \structa \apptruth \phi ( \vec a )$ (\structa  approximately satisfies
 $\phi$) 
iff 
\begin{equation*}
\exists h\in I(\phi(\vec x)) \mbox{ } \forall 
n\in\omega \mbox{, }
\structa \verify
([\phi(\vec a)]_{h})_{n} \findef
\end{equation*}
\end{defn}

{\bf Note\/} It is clear form the above definition that $\apptruth$ ``a la Henson'' and 
$\apptruth$ for $L_{A}$ coincide for formulas in $L_{PB}$. Hence, from now on, there shall 
be no confusion concerning the notion of $\apptruth$ being used.

\section{Uniformity Theorem for $L_{A}$}
\label{comply}

In this section we fix a countable signature $\Omega$.

We begin by proving that Henson's compactness theorem for 
$\apptruth$ in $L_{PB}$ in fact holds for $\apptruth$ in $L_{A}$.
This is 
not surprising, since the approximate formulas in $L_{A}$
 are positive bounded formulas.

Let us recall first three fundamental results
 for approximate truth in $L_{PB}$.
The interested reader can get details of the proofs
 in \cite{Henson&Iovino} or \cite{Iovino}.

\begin{thm} Henson's Compactness Theorem
\label{compact}

Let $\Sigma$ be a theory in $L_{PB}$, such that for every
 finite $F=\{\sigma_{i}: i \leq k\} 
\subseteq \Sigma$, for every integer $n$ there exists a
 normed space structure 
$E_{n}$ such that $E_{n} \models \bigwedge_{i=1}^{k}
 (\sigma_{i})_{n}$. Then there exists a 
normed space structure $E$ such that $E \apptruth \Sigma$.
\end{thm}

For the next theorem we need a definition.

\begin{defn} $\kappa$-saturated normed structures.

A normed space structure $E$ is $\kappa$-saturated if
 it approximately realizes any consistent set of
 formulas in $L_{PB}$ containing less than $\kappa$ constants
 and norm bounds for elements from
 $E$.\findef

\end{defn}

\begin{thm} $\aleph_{1}$-saturated structures.
\label{saturated}

For any normed structure $E$, there exists an
 approximate elementary extension $F$ of 
$E$ (i.e. $E$ and $F$ approximately satisfy
 the same formulas in $L_{PB}$ with parameters 
in $E$) that is $\aleph_{1}$-saturated.

\end{thm}

The final theorem shows that $\aleph_{1}$-saturated
 structures are "rich" for $L_{PB}$:

\begin{thm} 
\label{rich}

If $E$ is an $\aleph_{1}$-saturated and
 $\phi(\vec x) \in L_{PB}$ then $E \models 
\phi(\vec a)$ iff $E\apptruth \phi(\vec a)$.

\end{thm}

We use the above theorems to prove first that Henson's Compactness
 Theorem holds 
in fact for $\apptruth$ in $L_{A}$.

\begin{defn}
Let 
 $\Theta$ be a collection of sentences in $L_{A}$.
 We say that $\Theta$ is
 approximately finite consistent
iff there exists a set of branches $\Lambda=
\{h(\sigma) \in I(\sigma) : \sigma \in \Theta\}$ such that
for every finite subset $F \subset \Theta$, for every integer
$n$, there exists a normed space structure $\structn$ such that:
\begin{equation*}
\forall \sigma \in F, \mbox{ } \structn \verify ([\sigma]_{h(\sigma)})_{n}
\mbox{. \findef}
\end{equation*}
\end{defn}

Note that the above definition of approximate
finite consistency for $L_{A}$ restricted
to $L_{PB}$ coincide
with Henson's definition of finite consistency for $L_{PB}$
(\cite{Henson&Iovino}).

The main consequence of the previous lemmas is the following
Model Existence Theorem for $(L_{A}, \apptruth)$ and
$(L_{A},\models)$.

\begin{thm} 
\label{compactLA}
Fix 
 an approximately finite consistent collection $\Theta$ of
 sentences in $L_{A}$.
Then there exists a normed space structure $\structa $
 such that $E \apptruth \Theta$,
 and $E\models \Theta$. 
\end{thm}

\begin{proof}

Fix $\Theta$ as in the hypothesis of the theorem, and let 
 $\Lambda=\{h(\sigma) \in I(\sigma): \sigma \in \Theta\}$
 be the set of branches associated with 
$\Theta$. Let 
$
T=\{[\sigma]_{h(\sigma)} : \sigma \in \Theta\}
$. Clearly $T$ is a theory in $L_{PB}$ that satisfies
 the hypothesis of Henson's Compactness Theorem.
It follows from
Theorem~\ref{saturated}
that there exists an $\aleph_{1}$-saturated structure $E$ such that:
$ E \apptruth  T.$ By
the definition of $\apptruth$ in $L_{A}$ it follows that $E \apptruth
\Theta$. It remains to prove that 
$E \models \Theta$.

We claim:

For every formula $\phi (\vec x) \in L_{A}$, $E \apptruth \phi(\vec a) $
 iff 
$E \models \phi(\vec a)$.

proof: By induction in the formulas of $L_{A}$. The proofs 
for $L_{PB}$ and for the 
 countable (or finite) conjunction  steps
are direct, and are left to the reader.

{\bf Negation\/}.
$\Rightarrow$. Assume that $E \apptruth \neg \phi(\vec a)$.
 Assume also, in order 
to get a contradiction, that $E \models \phi(\vec a)$.
 By induction hypothesis it follows that 
$E \apptruth \phi(\vec a)$ which implies that there exists a
 branch $h \in I(\phi)$ such that 
\begin{equation}
E \models \bigwedge_{n=1}^{\infty} ([\phi]_{h})_{n}. \label{once}
\end{equation}

However, since $E \apptruth \neg \phi (\vec a)$,
 it follows from the definition of 
the approximate truth for the negation (Definition~\ref{defofapprox})
 that there exists a 
function $f=(f_{1},f_{2}): \omega \rightarrow I(\phi) \times \omega$
 with the following 
"weak" surjectivity property:
$$ \forall g \in I(\phi) \mbox{ } \exists s \mbox{, }
([\phi]_{g})_{f_{2}(s)}\equiv [(\phi]_{f_{1}(s)})_{f_{2}(s)}$$
and such that $E \models \bigwedge_{s=1}^{\infty}
 neg([\phi(\vec a)]_{f_{1}(s)},f_{2}(s))$.
From these two properties of $f$ it follows that
 there exists an $m$ such that 
$E \models neg ([\phi(\vec a)]_{h},m)$, which implies from
the properties of the weak approximate negation (Lemma~\ref{negando})
that
$E \not \models \bigwedge_{n=1}^{\infty} ([\phi]_{h})_{n}$, 
 but this contradicts the statement~\ref{once}.

$\Leftarrow$. Assume that $E \models \neg \phi (\vec a)$.
 By induction hypothesis, we 
get that $E \not \apptruth \phi(\vec a)$. 

Hence, for every branch $h \in I(\psi)$, there exists 
an integer $m$ such that $E \not \models
([\phi(\vec a)]_{h})_{m}$.
  We invoke now Lemma~\ref{negando} to obtain that there exists an
  integer $n$ such that 
 $E \models neg([\phi(\vec a)]_{h},n)$.

Consider now that the collection of all formulas of the form
$neg([\phi(\vec a)]_{h},n)$ that  hold in $E$. From Lemma~\ref{negando} it follows that those formulas are
 finitary  and belong to
 $L_{PB}$ so they are at most countable.
This implies that we can construct a function
 $f=(f_{1},f_{2}):\omega \rightarrow I(\phi) \times \omega$
with the "weak" surjectivity property (i.e. $\forall g \in I(\phi)$
$\exists s $  $([\phi]_{g})_{f_{2}(s)}\equiv
([\phi]_{f_{1}(s)})_{f_{2}(s)}$) and such that 
$E \models \bigwedge_{s=1}^{\infty}
 neg([\phi(\vec a)]_{f_{1}(s)},f_{2}(s))$.

We get then from the definition of approximate formulas for $L_{A}$
 (Definition~\ref{defofapprox}) that 
$E \apptruth \neg \phi(\vec a)$. This completes the proof of the
 negation step.

{\bf Existential \/}. There is only one interesting direction.
 Assume that 
$E \apptruth \exists \vec x (||\vec x|| \leq \vec r
 \wedge \phi(\vec a, \vec x))$. Then there 
exists $h \in I(\phi)$ such that
 $$E \models \bigwedge_{n=1}^{\infty} \exists \vec x(\bigwedge_{i=1}
^{Ind(n)} ||x_{i}|| \leq r_{i} 
\wedge ([\phi(\vec a, \vec x)]_{h})_{n}).$$
 Since the above formula is a countable conjunction of
finitary formulas in $L_{PB}$ and
 $E$ is $\aleph_{1}$-saturated,
 it follows that the conjunction 
$\bigwedge_{i=1}^{\infty} ||x_{i}|| \leq r_{i} \wedge
[ \phi(\vec a, \vec x)]_{h}$ is approximately realized in $E$ for some
$\vec b$.
This implies that $E \apptruth \phi(\vec a, \vec b)$, and hence,
 by induction hypothesis, $E \models \exists \vec x
(||\vec x|| \leq \vec {r} 
\wedge \phi(\vec a, \vec x))$.
 This completes the proof of the existential step and of the claim.

 From the above claim it follows that $E \models \Theta$. This completes
 the proof of the theorem.\end{proof}.

We are now ready to prove the main result of the paper:
the Uniformity Theorem. 

\begin{thm} Uniformity Theorem for $L_{A}$.

\label{uniformity}
Let $\Sigma$ be a theory in $L_{PB}$.
Let $\neg \phi$ be a sentence in $L_{A}$.
 If $\Sigma 
\models \neg \phi$ then for every $h \in I(\phi)$ there
 exists an integer $n$ such that 
$\Sigma \models \neg ([\phi]_{h})_{n}$.
\end{thm}

\begin{proof}
Assume, in order to get a contradiction, that there exists
 $h \in I(\phi)$ such that for every 
integer $n$ there is a normed structure $E_{n}$ satisfying:
$$ E_{n} \models \Sigma \mbox{ and } E_{n} \models ([\phi]_{h})_{n}$$

We can now invoke Theorem~\ref{natural} to obtain that there exists a
 normed structure 
$E$ such that $E\models \Sigma$ and $E\models \phi$,
 but this is a contradiction with the 
hypothesis.\end{proof}

\section{Applications}
\label{Applications}

A first corollary of the Uniformity Theorem concerns sentences of
the form $\sigma \Rightarrow \theta$ where $\sigma,\theta \in L_{PB}$.

\begin{cor} "Almost" Versions.
\label{almost}

Fix a signature $\Omega$. Suppose that $\Sigma$ is a theory in $L_{PB}$
such that $\Sigma \models (\sigma \Rightarrow \theta)$, for
$\sigma, \theta \in L_{PB}$. Then for every integer $n$ there
exists an integer $m$ such that:
$$ \Sigma \models \sigma_{m} \Rightarrow \theta_{n}$$
\end{cor}
\begin{proof}
It is enough to decode the approximate formulas corresponding to
$$\sigma \Rightarrow \theta : \neg (\sigma \wedge \neg \theta)$$

Since $\sigma, \theta \in L_{PB}$, then $I(\sigma)=\{\emptyset\}$
 and $I(\neg \theta)=(\{\emptyset\}\times\omega)^{\omega}$. It
follows that $I(\sigma \wedge \neg \theta)= \{\emptyset\}
\times (\{\emptyset\}\times \omega)^{\omega}$.

Note
that for every integer $n$,
the constant function $[n+1]: \omega \rightarrow \{\emptyset\}
\times \{n+1\}$ belongs to $I(\neg \theta)$. Hence it follows
 from the Uniformity Theorem (Theorem~\ref{uniformity}) that there
 exists an integer $m$ such that:
 $$ \Sigma \models \neg (\sigma_{m} \wedge
 ([\neg \theta]_{[n+1]})_{m})$$
 which decoded says
$$\Sigma \models \neg (\sigma_{m} \wedge
\bigwedge_{i=1}^{m} (neg(\theta,[n+1](i))\mbox{ })_{m} \mbox{ })$$
which is
$$\Sigma \models \sigma_{m} \Rightarrow
\neg \mbox{ }((neg(\theta,n+1))_{m}).$$ Invoking
now Lemma~\ref{negando} we get that 
$\Sigma \models \sigma_{m} \Rightarrow \theta_{n}$.\end{proof}

\begin{example} Ulam's Theorem
\label{Ulam}     

Note first that Ulam's Theorem is equivalent to the version
 where the target and
domain Banach space are the same.

 Fix now an arbitrary integer $k$ and let
$\Omega_{k}$ be the signature induced by any normed space structure
$\structa=(E, T)$ where $T:E \rightarrow E$ is a continuous
map with the property that $\forall x \in E$ $||T(x)|| \leq k ||x||$.

Consider, in $\Omega_{k}$, the theory
 $$\Sigma_{k}=\{\forall x (||x|| \leq n
 \Rightarrow ||T(x)|| \leq k||x||) :
 n \in \omega\}\wedge$$ $$ \{||T(0)||=0\} \wedge$$
$$ \{\forall x (||x||\leq n \mbox{ }
\Rightarrow \exists y (||y||\leq kn \wedge \mbox{ }||T(x)-y||=0)):
 n \in \omega\} $$
 that says that $T$  sends $0$ to $0$ and is an onto map. Clearly
 $\Sigma_{k} \subseteq L_{PB}$.

  Note that the version of Ulam's Theorem mentioned above
can be expressed as a formula in $L_{A}$:
 $$\phi: (\forall x,y (||x||,||y|| \leq 1 \Rightarrow ||T(x)-T(y)||=
||x-y||))$$
 $$ \Rightarrow (\forall x, y (||x||,||y|| \leq 1 \Rightarrow
 ||T(x)+T(y)-x-y||=0))$$
 which  is of the form $$\phi: (\sigma \Rightarrow \theta) :
 \neg (\sigma \wedge \neg \theta)$$
 with $\sigma,\theta \in L_{PB}$. Since $\Sigma_{k} \models
\phi$ we can invoke the
Corollary~\ref{almost}
  for formulas based on the signature $\Omega_{k}$ to
 obtain that for every integer $n$ there exists an integer $m$
  such that:

$$\Sigma_{k} \models [\forall x,y (||x||,||y|| \leq 1-1/m \Rightarrow$$
$$||x-y||-1/m \leq ||T(x)-T(y)|| \leq ||x-y||+1/m)]
\Rightarrow$$ $$ [\forall x,y (||x||,||y|| \leq 1-1/n \Rightarrow
||T(x)+T(y)-x-y|| \leq 1/n)]$$
which easily implies Gervitz's version.\findef
\end{example}

\begin{example} Behrends' Theorem.

Fix an integer $k$ and two rational numbers $1 \leq p,q \leq \infty$
and let $\Omega_{k}$ be the signature induced by
 a normed space structure $\structa=(E,P,Q, e_{1}, e_{2},f^{p},
f^{q})$ with:
\begin{itemize}
\item  $P,Q$ being 
linear projections (i.e. $P^{2}=P$ and $Q^{2}=Q$) with the property
that for every $x\in E$, $||P(x)|| \leq k||x||$ and $||Q(x)||\leq
k||x||$,
\item $e_{1},e_{2}$ being vectors with norm 1,
\item $f^{p}$ being a real valued function satisfying the same modulus
of continuity and the same bounds as the function $||.||^{p}$,
\item $f^{q}$ being a real valued function satisfying the same modulus
of continuity and the same bounds as the function $||.||^{q}$.
\end{itemize}

Consider, in $\Omega_{k}$, the theory $\Sigma_{k}=$
$$\{\forall x,y (||x||,||,y||\leq r \Rightarrow ||P(ax+by)-aP(x)-bP(y)||=0)
: a,b,r \in \mathbb{Q} \} \cup$$
$$\{\forall x,y (||x||,||y||\leq r \Rightarrow ||Q(ax+by)-aQ(x)-bQ(y)||=0)
: a,b,r \in \mathbb{Q} \} \cup$$
$$\{\forall x (||x||\leq r \Rightarrow ||P(P(x))-x||=0):
 r \in \mathbb{Q}\}\cup$$
$$\{\forall x (||x||\leq r \Rightarrow ||Q(Q(x))-x||=0):
 r \in \mathbb{Q}\}\cup$$
$$\{\forall x (||x||\leq r \Rightarrow (||x||\leq s \vee s^{p}
\leq f^{p}(x) \leq r^{p})) : s < r \mbox{ in } \mathbb{Q}\} \cup$$
$$\{\forall x (||x||\leq r \Rightarrow (||x||\leq s \vee s^{q}
\leq f^{q}(x) \leq r^{q})) : s < r \mbox{ in } \mathbb{Q}\} \cup$$
 $$\{ ||e_{1}||=1 \wedge ||e_{2}||=1 \wedge ||e_{1}-e_{2}|| \geq 1/2\}.$$
The last sentence listed for $\Sigma_{k}$ is equivalent
(using Riesz's Lemma)
to the statement that the dimension is $> 2$.
 Clearly, $\Sigma_{k} \subseteq L_{PB}$.

 Behrends' Theorem can be expressed by the following sentence in
 $L_{A}$:
 $$[\forall x (||x|| \leq 1 \Rightarrow (f^{p}(P(x))+f^{p}(x-P(x)))=
 f^{p}(x)) \wedge
\forall x (||x|| \leq 1$$
$$ \Rightarrow (f^{q}(P(x))+f^{q}(x-Q(x)))=
 f^{q}(x))] $$
$$\Rightarrow [ \forall x (||x||\leq 1 \Rightarrow ||P(Q(x)) -
Q(P(x))||=0 \wedge ||f^{p}(x)-f^{q}(x)||=0)]$$

It follows then from Corollary~\ref{almost} that for every integer $n$
there exists an integer $m$ such that for every normed space $E$
with $dim > 2$,: 
 $$[\forall x (||x|| \leq 1-1/m \Rightarrow $$ $$(f^{p}(x)-1/m \leq f^{p}(P(x))+
f^{p}(x-P(x)))\leq f^{p}(x)+1/m) \wedge$$
$$\forall x (||x|| \leq 1$$ $$ \Rightarrow (f^{q}(x) -1/m \leq
f^{q}(P(x))+f^{q}(x-Q(x)))\leq 
 f^{q}(x)+1/m)] $$
$$\Rightarrow [ \forall x (||x||\leq 1-1/n \Rightarrow ||P(Q(x)) -
Q(P(x))|| \leq 1/n\wedge ||f^{p}(x) -f^{q}(x)||\leq 1/n)].$$ 
from which one obtains the result from Cambern, Jaroz and Wodinski.\findef
\end{example}

The last example concerns a different type of formula:
 an infinitary one.

\begin{example} Krivine's Theorem.
\label{krivine}

A celebrated result by Krivine (\cite{Krivine}) states that for
every basic Schrauder sequence of unitary vectors
 $\{x_{n}\}_{n=1}^{\infty}$ in a normed space $E$,
there exists a $p \in [1,\infty)$ such that
 the usual basis of
 one of the spaces $\ell_{p}$  (or $c_{0}$)
 is block finitely representable in $\{x_{n}\}_{n=1}^{\infty}$.
We refer the reader to \cite{Lindenstrauss&Tzafriri} for the definition
and properties of the basic Schrauder sequences as well as
of the block finitely representable basis.

 Let $\Omega$ be the empty signature.  Let $\mathbb{Q}^{\#}=
 (\mathbb{Q}\cap [1,\infty])$. For a fix integer $n$, let
 $\vec {q}$ denote a vector $(q_{1},q_{2},\ldots,q_{n+1})$ of
 integers such that
 $q_{1} < q_{2} < \ldots < q_{n+1}$ and let $V_{n} \subset \omega^{n+1}$
 be the collection of all such vectors.
 
A weaker version of this theorem has the form: For every $K \geq 1$,
for every $\epsilon > 0$, for every $n \in \omega$, for every
normed space $E$,
$$E \models \forall \vec x (||\vec x|| \leq 1 \Rightarrow $$
$$[BaseK(\vec x) \Rightarrow \bigvee_{p \in Q^{\#}}
\bigvee_{\vec q \in V_{n}} \bigvee_{\vec b \in Q^{q_{n}+1}}
 \theta^{n,p,\epsilon}(\sum_{i=q_{1}+1}^{q_{2}} b_{i}x_{i},
\sum_{i=q_{2}+1}^{q_{3}} b_{i}x_{i},\ldots. \sum_{i=q_{n}+1}^{q_{n+1}}
b_{i}x_{i})]$$
where
\begin{itemize}
\item
$$BaseK(\vec x):\bigwedge_{i=1}^{\infty}
||x_{i}|| =1 \wedge \bigwedge_{n,m \in \omega} \bigwedge_{\vec a\in Q^{n+m}}
||\sum_{i=1}^{n} a_{i}x_{i}|| \leq K ||\sum_{i=1}^{n+m} a_{i}x_{i}||$$
is a positive formula that
states that the sequence $\{x_{i}\}_{i=1}^{\infty}$ is a basic
Schrauder sequence with constant $K$;
\item  for every integers
$n$ and $p \in Q^{\#}$,
$$ \theta^{n,p,\epsilon}(y_{1},\ldots, y_{n}): \bigwedge_{\vec c \in Q^{n}}
\sqrt[p]{\sum_{i=1}^{n} |c_{i}|^{p}} \leq ||\sum_{j=1}^{n} c_{j}y_{j}||
\leq (1+\epsilon)\sqrt[p]{\sum_{i=1}^{n} |c_{i}|^{p}}$$
states that the usual basis of the space $\ell_{p}^{n}$ is
$(1+\epsilon)$-equivalent to the vectors $(y_{1},\ldots,y_{n})$.
\end{itemize}

Note that the weaker form of Krivine's Theorem can be written as:
For every $K \geq 1$, for every $\epsilon > 0$,
for every integer $n$, 
$$ \models \neg \exists \vec x (||\vec x|| \leq 1 \wedge $$
$$[BaseK(\vec x) \wedge \bigwedge_{p \in \mathbb{Q}^{\#}}
\bigwedge_{\vec q \in V_{n}} \bigwedge_{\vec b \in Q^{q_{n}+1}}
 \neg \theta^{n,p,\epsilon}(\sum_{i=q_{1}+1}^{q_{2}} b_{i}x_{i},
\sum_{i=q_{2}+1}^{q_{3}} b_{i}x_{i},\ldots. \sum_{i=q_{n}+1}^{q_{n+1}}
b_{i}x_{i})]$$

Note also that, by virtue of the finite dimensionality of the
$\ell_{p}^{n}$, for every integer $n$, for every $p \in \mathbb{Q}^{\#}$,
for every $\epsilon >0$ there exists an integer $w(n,p,\epsilon)$
such that
$$\models (\theta^{n,p,\epsilon})_{w(n,p,\epsilon)} \Rightarrow
 \theta^{n,p, 2\epsilon}$$

Using the above remark, define for any $\epsilon \geq 0$, for
 every integer $n$, for every
$K \geq 1$, the function
$$h:\mathbb{Q}^{\#} \times V_{n}\times Q^{q_{n+1}} \rightarrow
\omega$$ such that $h(p,(q_{1},q_{2},\ldots,q_{n+1}),\vec b)=
w(n,p,\epsilon)+ 1$. We leave to the reader the verification that the pair
 $(\emptyset,h)$ is a branch of the tree of approximations of
 the formula
 $$BaseK(\vec x) \wedge
 \bigwedge_{p \in \mathbb{Q}^{\#}}
\bigwedge_{\vec q \in V_{n}} \bigwedge_{\vec b \in Q^{q_{n}+1}}
 \neg \theta^{n,p,\epsilon}(\sum_{i=q_{1}+1}^{q_{2}} b_{i}x_{i},
\sum_{i=q_{2}+1}^{q_{3}} b_{i}x_{i},\ldots. \sum_{i=q_{n}+1}^{q_{n+1}}
b_{i}x_{i})$$

For every $K\geq 1$, for every $\epsilon >0$, for every integer $n$
we can invoke the Uniformity Theorem for the branch
$(\emptyset,h)$ to obtain that there exists an integer $r$ such that:
$$
\models \neg \exists \vec x (||\vec x|| \leq 1+(1/r) \wedge
[(BaseK(\vec x))_{r} \wedge 
\bigwedge_{p \in \mathbb{Q}^{\#}\uparrow r}
\bigwedge_{\vec q \in V_{n}\uparrow r} \bigwedge_{\vec b \in Q^{q_{n}+1}
\uparrow r}$$ $$
 neg(\theta^{n,p,\epsilon}
(\sum_{i=q_{1}+1}^{q_{2}} b_{i}x_{i},
\sum_{i=q_{2}+1}^{q_{3}} b_{i}x_{i},\ldots. \sum_{i=q_{n}+1}^{q_{n+1}}
b_{i}x_{i}),h(p,\vec q,\vec b))]$$
which implies, using again Lemma~\ref{negando}, the property of $w(n,p,\epsilon)$ and the definition of $h$,
that 
$$ \models \forall \vec x (||\vec x|| \leq 1 \Rightarrow 
[(BaseK(\vec x))_{r} \Rightarrow $$ $$\bigvee_{p \in \mathbb{Q}^{\#}\uparrow r}
\bigvee_{\vec q \in V_{n}\uparrow r} \bigvee_{\vec b \in Q^{q_{n}+1}
\uparrow r}
 \theta^{n,p,2\epsilon}(\sum_{i=q_{1}+1}^{q_{2}} b_{i}x_{i},
\sum_{i=q_{2}+1}^{q_{3}} b_{i}x_{i},\ldots. \sum_{i=q_{n}+1}^{q_{n+1}}
b_{i}x_{i})].$$

This last statement can be written as follows:

For every $K > 1$, for every $\epsilon > 0$, for every integer $n$,
there exists a finite collection $I=\{p_{1},\ldots,p_{r}\}
\in Q\cap [1,\infty]$ and an integer
$m$ such that
for every finite basic sequence
$(x_{i})_{i=1}^{m}$ with basic constant $K$ in any normed space,
there exists a
$p \in I$ and a block basic sequence that is $1+\epsilon$
 equivalent to the usual basis of $\ell_{p}^{n}$.

 Compare this result with the Uniform Version of Krivine's
 Theorem obtained by Rosenthal (\cite{Rosenthal}):

 Fix arbitrary $K \geq 1$, $n \in \omega$ and $\epsilon >0$.
 There exists an $m$ such that if $(x_{i})_{i=1}^{m}$ is a finite
 basic sequence
  in any Banach space with basis constant $K$, then there exists
  $1 \leq p \leq \infty$ and a block sequence $(y_{i})_{i=1}^{n}$
  so that $(y_{i})_{i=1}^{n}$ is $(1+\epsilon)$-isomorphic
 to the unit vector
  basis of $\ell_{p}^{n}$.\findef
  \end{example}

\bigskip

Beaver College, Glenside, PA.

ortiz@beaver.edu
\medskip

Mathematical Subject Classification, 1991: 

 Primary 03C65

 Secondary  46B08, 46B20
\medskip

Keywords: approximate truth, compactness theorem, infinitary logic, 
normed space structures, uniformity results.


\begin{thebibliography}{99}


\frenchspacing
\parskip 12pt



\bibitem{Behrends} E. Behrends, {\it $L^{P}$-struktur in
Banachr\"{a}umen\/}, {\bf Studia Math.\/}, vol 55, (1976), pp. 71-85.

\bibitem{Cambern} M. Cambern, K. Jarosz and G. Wodinski,
{\it Almost $L^{p}$-projections and $L^{p}$-isomorphisms\/},
{\bf Proc. Royal. Soc. Edimburgh\/}, vol 113A, (1989), pp. 13-25.




\bibitem{Gevirtz} J. Gevirtz, {\it Surjectivity in Banach spaces and the
 Mazur-Ulam theorem on isometries\/}, {\bf Trans. Amer. Math. Soc.\/},
 vol 274, (1982), pp. 307-318.
 



\bibitem{Henson} C.W. Henson, {\it When Do Two Banach Spaces Have
 Isometrically
Isomorphic
 Nonstandard Hulls?\/}, {\bf Israel J. Math.}, vol 22, (1975), pp. 57-67.




   \bibitem{Henson&Iovino} C.W. Henson and J. Iovino, {\it
Banach Space Model
   Theory, I:
   Basics\/}, in preparation.

   \bibitem{Hodges} W. Hodges, {\bf Model Theory.\/} Cambridge
   University Press, Cambridge, 1994.
   
\bibitem{Iovino} J. Iovino,
{\bf Stable Theories in Functional Analysis\/}, 
Ph.D. dissertation, University of Illinois at Urbana-Champaign, 1994.

\bibitem{James} R.C. James, {\it Weak Compactness and Reflexivity\/},
 {\bf
Israel J. Math.\/}, vol. 2, (1964), pp. 101-119.

\bibitem{Jarosz} K. Jarosz, {\it Ultraproducts and Small
 Bound Perturbations\/}, {\bf Pacific Journal of Mathematics\/},
 vol 148, (1991), pp. 81-88.
 
\bibitem{Krivine} J.L.Krivine, {\it Sous espaces de dimension finie des 
espaces de Banach reticul\'{e}s\/}, {\bf Ann. of Math.\/},
vol 104, 
(1974), pp. 213-253.

 \bibitem{Lindenstrauss&Tzafriri} J. Lindenstrauss and
L. Tzafriri, {\bf
   Classical Banach Spaces\/}, Springer Verlag, Berlin, 1977.
 
\bibitem{Ortiz} C.E. Ortiz, {\it Approximate Truth and
 Nonstandard Analysis\/},
{\bf Cahiers du Centre de Logique\/}, vol. 9, Academia Bruylant,
 Louvain la Neuve, 1996.

 \bibitem{Rosenthal} H. Rosenthal, {\it On a theorem of Krivine concerning 
block finite representability of $\ell_{p}$ in general Banach spaces\/},
 {\bf Journal of Functional Analysis\/}, vol 28, 
(1978), pp. 197-225.

  \bibitem{Van Dulst} D. Van Dulst, 
  {\bf Reflexive and Superreflexive Banach Spaces\/},
  Mathematical Centre Tract, Amsterdam, 1978.


 \end{thebibliography}
\end{document}